\documentclass{amsart}
\usepackage{amsmath,amscd,amssymb,amsthm,graphicx}

\def\cal{\mathcal}

\def\C{{\cal C}} 
\def\D{{\cal D}} 
\def\E{{\cal E}}

\def\I{{\cal I}}

\def\L{{\cal L}}

\def\S{{\cal S}}

\def\frak{\mathfrak}

\def\da{{\frak a}}

\def\dg{{\frak g}}
\def\dgl{\dg\dl}

\def\dk{{\frak k}}
\def\dl{{\frak l}}
\def\dm{{\frak m}}

\def\dgo{{\frak o}}

\def\dgp{{\frak p}}

\def\dr{{\frak r}}
\def\ds{{\frak s}}
\def\dsl{\ds\dl}
\def\dsp{\ds\dgp}

\def\dU{{\frak U}}

\def\dz{{\frak z}}

\def\dZ{{\frak Z}}

\def\Bbb{\mathbb}

\def\bC{\Bbb C}
\def\bD{\Bbb D}

\def\bH{\Bbb H}

\def\bN{\Bbb N}

\def\bR{\Bbb R}

\def\bZ{\Bbb Z}

\def\ep{\epsilon}

\def\ot{\otimes}

\def\t{\tilde}

\def\ad{\mathop{\rm ad}\nolimits}

\def\and{\mbox{\rm \ and\ }}

\def\Aut{\mathop{\rm Aut}\nolimits}

\def\det{\mathop{\rm det}\nolimits}

\def\oh{{\ts\frac{1}{2}}}

\def\scrm{\scriptsize\rm}
\def\SL{\mathop{\rm SL}\nolimits}

\def\SO{\mathop{\rm SO}\nolimits}

\def\Sp{\mathop{\rm Sp}\nolimits}
\def\Span{\mathop{\rm Span}\nolimits}

\def\SU{\mathop{\rm SU}\nolimits}

\def\thup{{\mbox{\scrm th}}}
\def\tr{\mathop{\rm tr}\nolimits}

\def\ts{\textstyle}

\def\Vec{\mathop{\rm Vect}\nolimits}

\def\dog{differential operator}

\def\ido{invariant \dog}
\def\ie{{\em i.e.,\/}}
\def\iff{if and only if}

\def\irr{irreducible}
\def\irrep{irreducible representation}

\def\meno{\medbreak\noindent}

\def\r{representation}

\def\sq{subquotient}

\def\uea{universal enveloping algebra}

\def\vb{vector bundle}

\def\vs1{$\Vec(S^1)$}

\def\iar{irreducible admissible representation}
\def\Lie{\mathop{\rm Lie}\nolimits}

\def\dzi{\det(Z)^{-1}}
\def\zih{Z^{-1/2}}

\newtheorem{lemma}{Lemma}[section]

\newtheorem{prop}[lemma]{Proposition}
\newtheorem{thm}[lemma]{Theorem}
\newtheorem{cor}[lemma]{Corollary}

\title[Centers and characters of IDO algebras]{Centers and characters of Jacobi group-invariant differential operator algebras}

\author{Charles H.\ Conley}
\address{Department of Mathematics \\University of North Texas \\Denton TX 76203 \\USA} 
\email{conley@unt.edu}

\author{Rabin Dahal}
\address{Department of Mathematics \\Texas College \\Tyler TX 75702 \\USA} 
\email{rdahal@texascollege.edu}
\thanks{The first author was partially supported by Simons Foundation Collaboration Grant 207736.}

\begin{document}

\begin{abstract}
We study the algebras of \dog s invariant with respect to the scalar slash actions of real Jacobi groups of arbitrary rank.  These algebras are non-commutative and are generated by their elements of orders~2 and~3.  We prove that their centers are polynomial in one variable and are generated by the Casimir operator.  For slash actions with invertible indices we also compute the characters of the IDO algebras: in rank exceeding~1 there are two, and in rank~1 there are in general five.  In rank~1 we compute in addition all \iar s of the IDO algebras.
\end{abstract}


\maketitle

% \noindent Mathematics Subject Classifications: 22E47, 17B10, 11F50

\section{Introduction}  \label{Intro}

The real Jacobi group $G^J_N$ of degree~$1$ and rank $N$ is the semidirect product of $\SL_2(\bR)$ acting on a certain Heisenberg central extension of $N$ copies of the standard \r.  The theory of Jacobi forms on $G^J_1$ was developed in \cite{EZ85} and has since been widely used.  The generalization to Jacobi forms of higher degree and rank, and in particular to Jacobi forms on $G^J_N$, was initiated in \cite{Zi89}.

Jacobi forms are holomorphic by definition.  In the definition of automorphic forms on general semisimple Lie groups, holomorphicity is replaced by the condition that the form be an eigenfunction of all elements of the center of the \uea\ of the group \cite{HC59}, or at least that it generate a finite dimensional \r\ of the center \cite{Bo66}.  Berndt and Schmidt \cite{BS98} proposed a similar approach for the Jacobi group, which of course is not semisimple: they defined automorphic forms on $G^J_1$ using the cubic invariant operator $C$ defined in \cite{BB90} in place of the center.  This definition was used in the study of non-holomorphic Maa\ss-Jacobi forms in \cite{Pi09} and \cite{BR10}.  At this time it was realized that the center of the \uea\ of $G^J_1$ in fact acts by $\bC[C]$ \cite{BCR12}.  Thus the definition from \cite{BS98} coincides with the classical one.

In order to describe the \r-theoretic framework of automorphic forms, let $G$ be a connected real Lie group, $K$ a Lie subgroup, and $V$ a complex \r\ of $K$.  In this setting we have the $G$-\vb\ $G \times_K V$ over the homogeneous space $G/K$, and the space $C^\infty_{\sec}(G \times_K V)$ of smooth sections of $G \times_K V$.  The group $G$ acts on $C^\infty_{\sec}(G \times_K V)$ by translations, its complexified Lie algebra $\dg$ acts by vector fields, and its \uea\ $\dU(\dg)$ acts by \dog s.  In the event that the \vb\ is topologically trivial, its sections may be identified with the $V$-valued functions.  In this case the associated right action of $G$ on the functions is called a {\em slash action.\/}

Observe that the center $\dZ(\dg)$ of $\dU(\dg)$ acts on $C^\infty_{\sec}(G \times_K V)$ by \dog s which commute with the action of $G$.  In general, such differential operators are said to be {\em invariant.\/}  The {\em invariant \dog\ algebra\/} of $G \times_K V$ (its ``IDO algebra'') is the associative algebra of all invariant \dog s on $C^\infty_{\sec}(G \times_K V)$.  We will denote it by
\begin{equation*}
   \bD(G \times_K V).
\end{equation*}

This algebra is not necessarily commutative, and we write $Z \bigl(\bD(G \times_K V)\bigr)$ for its center.  It is a crucial point that the $\dZ(\dg)$-action manifestly commutes with $\bD(G \times_K V)$ as well as with the $G$-action, so there is a natural homomorphism
\begin{equation*}
   \dZ(\dg) \to Z \bigl(\bD(G \times_K V)\bigr).
\end{equation*}
In general, this homomorphism is neither injective nor surjective.

IDO algebras have been the subject of many deep investigations; see for example the survey articles \cite{He77} and \cite{He79} and the references therein.  One focus is on {\it eigenspace \r s:\/} if $\chi: \bD(G \times_K V) \to \bC$ is a {\em character,\/} \ie\ a homomorphism onto the scalars, the associated eigenspace \r\ of $G$ is
\begin{equation*}
   C^\infty_{\sec}(G \times_K V)_\chi := \bigl\{s \in C^\infty_{\sec}(G \times_K V):
   Ds = \chi(D) s\ \forall\ D \in \bD(G \times_K V) \bigr\}.
\end{equation*}

In the case of primary importance, $G$ is a connected non-compact semisimple group with finite center, $K$ is a maximal compact subgroup, and $V$ is the trivial \r\ $\bC$.  Here it is known that $\bD(G \times_K \bC)$ is commutative (see \cite{He79} for the origin of this result, as well as an extension to more general symmetric spaces).  Moreover, in all but four cases $\bD(G \times_K \bC)$ is in fact the image of $\dZ(\dg)$, and even in those four cases all elements of $\bD(G \times_K \bC)$ can be expressed as ratios of images of elements of $\dZ(\dg)$ \cite{He92}.  In the event that $K$ has an abelian factor and hence non-trivial 1-dimensional \r s, $\bD(G \times_K V)$ is still commutative for $V$ 1-dimensional, because the homomorphism from $\dZ(\dg)$ to $\bD(G \times_K V)$ preserves degree and the graded commutative symbol algebras of $\bD(G \times_K \bC)$ and $\bD(G \times_K V)$ are the same.  However, $\bD(G \times_K V)$ is in general not commutative when $V$ is multidimensional.  In the sense of \cite{HC59} and \cite{Bo66}, an automorphic form on $G$ is a section of $G \times_K V$ that is invariant under a discrete subgroup of $G$, is an eigenfunction of all elements of $\dZ(\dg)$, and satisfies a growth condition.

For general $G$, speaking heuristically, if $\bD(G \times_K V)$ contains ``enough'' of the commutant of the action of $G$, one might expect that $C^\infty_{\sec}(G \times_K V)$ decomposes into a multiplicity-free direct sum (possibly partially continuous) of \irrep s under the joint action of $G$ and $\bD(G \times_K V)$, \ie\ under $\dU(\dg) \ot \bD(G \times_K V)$.  In particular, one may expect the eigenspace \r s of $G$ to be \irr.  As described in \cite{He77}, in many settings this is true for most but not all characters of $\bD(G \times_K V)$.

Jacobi forms of degree~$1$ and rank~$N$ live on the homogeneous space $G^J_N/K^J_N$, where $K^J_N$ is the product of $\SO_2$ with the center $Z(G^J_N)$.  In this article we investigate the algebraic structure of $\bD(G^J_N \times_{K^J_N} V)$.  Since $K^J_N$ is abelian, we will assume that $V$ is 1-dimensional.  Moreover, we will only treat those $V$ with {\em invertible indices\/} (see Section~\ref{D}), the situation of interest in number theory.

In contrast with the semisimple case, $\bD(G^J_N \times_{K^J_N} V)$ is not commutative (this goes back to Berndt: see \cite{BB90}).  Our first result is that for all $N$, the center $Z\bigl(\bD(G^J_N \times_{K^J_N} V)\bigr)$ is precisely the image of $\dZ(\dg^J_N)$.  This was previously known for $N=1$ \cite{BCR12} and $N=2$ \cite{Da13}.  The structure of $\dZ(\dg^J_N)$ was deduced for $N=1$ in \cite{BCR12} and for $N > 1$ in \cite{CR}: it is generated by the center $\dz(\dg^J_N)$ of $\dg^J_N$ itself and in addition a single ``Casimir element'' of degree $N+2$, which acts by an IDO of degree~$3$ for $N=1$ and of degree~$4$ for $N > 1$.

Our second result is the description of the characters of $\bD(G^J_N \times_{K^J_N} V)$: for $N > 1$ there are two, and for $N = 1$ there are five for most slash actions, but only four for certain special slash actions.  For $N = 1$ we also classify the \irr\ {\em admissible \r s\/} of $\bD(G^J_1 \times_{K^J_1} V)$ (see Section~\ref{Reps}), which include its \irr\ finite dimensional \r s.  This is accomplished by identifying $\bD(G^J_1 \times_{K^J_1} V)$ with a subalgebra of $\dU(\dsl_2)$.  It could be interesting to determine which of these characters and higher dimensional \r s occur as constituents of the action of $\bD(G^J_N \times_{K^J_N} V)$ on $C^\infty_{\sec}(G^J_N \times_{K^J_N} V)$.

We begin in Section~\ref{GJN} by recalling the structure of the Jacobi Lie group.  In Section~\ref{D} we define the relevant homogeneous space, its scalar slash actions, and their IDO algebras.  Most of these first two sections is a recapitulation of material in \cite{BCR12} and \cite{CR}.  Our results are given in Sections~\ref{Z}, \ref{Characters}, and~\ref{Reps}.  We conclude in Section~\ref{Remarks} with some remarks on directions for further research.

\section{The Jacobi Lie group $G^J_N$}  \label{GJN}

Throughout this article all vector spaces are complex unless stated otherwise.  We will write $\bN$ for the non-negative integers and $\bZ^+$ for the positive integers.

Let $M_{m,n}(R)$ be the set of $m \times n$ matrices over a ring $R$, abbreviate $M_{n,n}(R)$ by $M_n(R)$, and write $M^T_n(R)$ for the subset of symmetric matrices.  Given a matrix $A$, write $A^T$ for its transpose, and if it is square, $\tr(A)$ and $\det(A)$ for its trace and determinant.  Let $I_n$ be the identity matrix in $M_n(R)$, and set $J := \bigl(\begin{smallmatrix} 0 & -1 \\ 1 & \ \ 0\end{smallmatrix} \bigr)$. 

The {\em real Jacobi group\/} $G^J_N$ of rank~$N$ is defined by
\begin{equation*}
   G^J_N := \bigl\{ (M, X, \kappa): M \in \SL_2(\bR),\ X \in M_{N,2}(\bR),\ \kappa + \oh X J X^T \in M_N^T(\bR) \bigr\},
\end{equation*}
\begin{equation*}
   (M_1, X_1, \kappa_1) (M_2, X_2, \kappa_2) :=
   (M_1 M_2, X_1 M_2 + X_2, \kappa_1 + \kappa_2 - X_1 M_2 J X_2^T).
\end{equation*}

\subsection{The Jacobi Lie algebra $\dg^J_N$} \label{Sec 2.1}

The complexified Lie algebra of $G^J_N$ is
\begin{equation*}
   \dg^J_N = \bigl\{ (M, X, \kappa): M \in \dsl_2(\bC),\ X \in M_{N,2}(\bC),\ \kappa \in M_N^T(\bC) \bigr\},
\end{equation*}
\begin{equation*}
   \bigl[ (M_1, X_1, \kappa_1), (M_2, X_2, \kappa_2) \bigr] =
   \bigl( [M_1, M_2], X_1 M_2 - X_2 M_1, X_2 J X_1^T - X_1 J X_2^T \bigr).
\end{equation*}

In Section~5 of \cite{CR} one may find the formula for the exponential map of $G^J_N$, as well as an embedding of $G^J_N$ in a minimal parabolic subgroup of $\Sp_{2N+2}(\bR)$.  (We remark that in that paper, the Jacobi group is called $\t G^J_N$, and the symbol $G^J_N$ is used for its centerless quotient.)  The centers of $G^J_N$ and $\dg^J_N$ are
\begin{equation*}
   Z(G^J_N) = \{I_2\} \times \{0\} \times M^T_N(\bR), \qquad
   \dz(\dg^J_N) = \{0\} \times \{0\} \times M^T_N(\bC).
\end{equation*}

We now fix a basis for $\dg^J_N$.  Let $\ep_r$ be the $r^\thup$ standard basis vector of $R^m$, and let $\ep_{rs}$ be the elementary matrix of $M_{m,n}(R)$ with $(r,s)$ entry~1 and other entries~0.  Set $E := \ep_{12}$, $F := \ep_{21}$, and $H := \ep_{11} - \ep_{22}$, the standard basis of $\dsl_2$.  For a basis of $M_{N,2}$, take $e_r := \ep_{r,2}$ and $f_r := \ep_{r,1}$, and for a basis of $M_N^T$, take $Z_{rs} := \oh (\ep_{rs} + \ep_{sr})$.  Then $\dg^J_N$ has basis
\begin{equation} \label{gJN basis}
   \bigl\{ E, F, H;\ e_r, f_r,\, 1 \le r \le N;\ Z_{rs},\, 1 \le r \le s \le N \bigr\}.
\end{equation}
The brackets of this basis are as follows: those on $\dsl_2$ are standard, $[e_r, f_s] = -2 Z_{rs}$, and as noted, $M^T_N$ is central.  Under $\ad(H)$, the $e_r$ are of weight~$1$ and the $f_r$ are of weight~$-1$.  The $\ad(E)$ and $\ad(F)$ actions are given by
\begin{equation*}
   \ad(E):\ e_r \mapsto 0,\ f_r \mapsto -e_r; \qquad
   \ad(F):\ e_r \mapsto -f_r,\ f_r \mapsto 0.
\end{equation*}

There is an automorphism $\theta$ of $\dg^J_N$ of order~4, defined by
\begin{equation*}
  \theta:\ H \mapsto -H, \quad E \mapsto -F \mapsto E, \quad
  e_r \mapsto f_r \mapsto -e_r, \quad Z_{rs} \mapsto Z_{rs}.
\end{equation*}

\subsection{The center $\dZ(\dg^J_N)$ of $\dU(\dg^J_N)$}

It will be useful to write $e$ and $f$ for the column vectors with entries $e_r$ and $f_r$, respectively, and $Z$ for the symmetric matrix with entries $Z_{rs}$.  Thus for example $f^T Z e$ and $\det(Z)$ are elements of the \uea\ $\dU(\dg^J_N)$.

Consider the localization $\dU(\dg^J_N)[\dzi]$ of $\dU(\dg^J_N)$ at $\det(Z)$.  Observe that it contains the entries of $Z^{-1}$, and define a map $\nu: \dsl_2 \to \dU(\dg^J_N)[\dzi]$ by
\begin{eqnarray*}
   \nu(H) &:=& H - {\ts\frac{1}{4}} (f^T Z^{-1} e + e^T Z^{-1} f)\ =\ H - \oh (f^T Z^{-1} e - N), \\[6pt]
   \nu(E) &:=& E - {\ts\frac{1}{4}} e^T Z^{-1} e, \qquad
   \nu(F)\ :=\ F + {\ts\frac{1}{4}} f^T Z^{-1} f.
\end{eqnarray*}
The following factorization of $\dU(\dg^J_N)[\dzi]$ is given in Section~5 of \cite{CR}.  Note that the radical of $\dg^J_N$ is the Heisenberg Lie algebra
\begin{equation*}
   \dr^J_N := \Span \bigl\{ e_r, f_r,\, 1 \le r \le N;\ Z_{rs},\, 1 \le r \le s \le N \bigr\}.
\end{equation*}

\begin{thm} \label{factorization} \begin{enumerate}
\item[(i)] $\nu$ is a Lie algebra homomorphism.  It extends to an associative algebra homomorphism $\nu: \dU(\dsl_2) \to \dU(\dg^J_N)[\dzi]$.
\item[(ii)] The image $\nu\bigl(\dU(\dsl_2)\bigr)$ commutes with $\dU(\dr^J_N)$.
\item[(iii)] The natural multiplication map defines an associative algebra isomorphism
\begin{equation*}
   \nu\bigl(\dU(\dsl_2)\bigr) \otimes \dU(\dr^J_N)[\dzi] \cong \dU(\dg^J_N)[\dzi].
\end{equation*}
\end{enumerate}
\end{thm}

This factorization is used in \cite{CR} to compute the center $\dZ(\dg^J_N)$ of $\dU(\dg^J_N)$; we give the result in the next proposition.  Recall that the Casimir element of $\dU(\dsl_2)$ is
\begin{equation*}
   \Omega_{\dsl_2} := H^2 + 2H + 4FE = H^2 - 2H + 4EF.
\end{equation*}
Since $\nu(\Omega_{\dsl_2})$ is quadratic in $Z^{-1}$, multiplying by $\det(Z)^2$ clears its denominator: $\det(Z)^2 \nu(\Omega_{\dsl_2})$ is in $\dU(\dg^J_N)$.  In fact, as explained in \cite{CR}, its denominator is $\det(Z)$, \ie\ $\det(Z) \nu(\Omega_{\dsl_2})$ is in $\dU(\dg^J_N)$.  In order to maintain the normalization of \cite{CR}, we define the {\em Casimir element of\/} $\dU(\dg^J_N)$ to be
\begin{equation*}
   \Omega_N := \det(Z) \bigl( \nu(\Omega_{\dsl_2}) - {\ts\frac{1}{4}} N (N + 4) \bigr).
\end{equation*}

\begin{prop} \label{center}
\begin{enumerate}
\item[(i)]  $\dZ(\dg^J_N) = \dU\bigl(\dz(\dg^J_N)\bigr) [\Omega_N] = \bC[\Omega_N, Z_{rs}: 1 \le r \le s \le N]$.
\item[(ii)]  $\Omega_N$ is in the augmentation ideal of\/ $\dU(\dg^J_N)$.
\end{enumerate}
\end{prop}

\subsection{The structure of $\dU(\dr^J_N)$} \label{UrJN}

Suppose now that $\dU(\dg^J_N)[\zih]$ is a further extension of $\dU(\dg^J_N)$ in which $Z$ has both a reciprocal and a symmetric square root.  Let us observe that the symplectic Lie algebra $\dsp_{2N}$ embeds in $\dU(\dr^J_N)[\zih]$.  To explain, recall that
\begin{equation*}
   \dsp_{2N} := \Bigl\{ \Bigl(\begin{matrix} A & \ B \\ C & -A^T \end{matrix} \Bigr):
   A \in M_N(\bC);\, B, C \in M_N^T(\bC) \Bigr\}.
\end{equation*}
For a basis of $\dsp_{2N}$, take 
\begin{equation*}
   A_{rs} := \Bigl(\begin{matrix} \ep_{rs} & \ 0 \\ 0 & -\ep_{sr} \end{matrix} \Bigr), \quad
   B_{rs} := \Bigl(\begin{matrix} 0 & \oh(\ep_{rs} + \ep_{sr}) \\ 0 & 0 \end{matrix} \Bigr), \quad
   C_{rs} := B_{rs}^T.
\end{equation*}
Note that $f e^T - (e f^T)^T = 2Z$.  The embedding $\dsp_{2N} \hookrightarrow \dU(\dr^J_N)[\zih]$ is
\begin{equation*}
   A_{rs} \mapsto -{\ts\frac{1}{4}}
   \bigl( (\zih f e^T \zih)_{rs} + (\zih e f^T \zih)_{sr} \bigr)
   = \oh \bigl( I_N - \zih f e^T \zih \bigr)_{rs},
\end{equation*}
\begin{equation*}
   B_{rs} \mapsto -{\ts\frac{1}{4}} (\zih f f^T \zih)_{rs}, \quad
   C_{rs} \mapsto {\ts\frac{1}{4}} (\zih e e^T \zih)_{rs}.
\end{equation*}

\section{The homogeneous space $G^J_N / K^J_N$ and its slash actions} \label{D}

The homogeneous space of $G^J_N$ of interest in this article is $\bH \times \bC^N$, where $\bH$ is the Poincar\'e upper half plane.  The explicit action may be found in \cite{CR}; in this paper we will only need the associated stabilizer subgroup of $G^J_N$, which is
\begin{equation*}
   K^J_N := \bigl\{ (M, 0, \kappa): M \in \SO_2,\ \kappa \in M_N^T(\bR) \bigr\}.
\end{equation*}
(Like $G^J_N$, $K^J_N$ is called $\t K^J_N$ in \cite{CR}, the symbol $K^J_N$ being used for the quotient of the stabilizer subgroup by the center of the Jacobi group.)

The complexified Lie algebra of $K^J_N$ is
\begin{equation*}
   \dk^J_N = \bigl\{ (M, 0, \kappa): M \in \dgo_2(\bC),\ \kappa \in M_N^T(\bC) \bigr\}.
\end{equation*}
The basis~(\ref{gJN basis}) is not compatible with $\dk^J_N$, so we recall from \cite{CR} the {\em tilde basis\/} of $\dg^J_N$:
\begin{equation*}  \begin{array}{rclrclrcl}
   \t H &:=& i(F-E), \ &
   \t E &:=& {\ts\frac{1}{2}} \bigl(H + i(F + E)\bigr), \ &
   \t F &:=& {\ts\frac{1}{2}} \bigl(H - i(F + E)\bigr), \\[6pt]
   \t Z_{rs} &:=& {\ts\frac{1}{2}} iZ_{rs}, \ &
   \t e_r &:=& {\ts\frac{1}{2}} (f_r + ie_r), \ &
   \t f_r &:=& {\ts\frac{1}{2}} (f_r - ie_r).
\end{array} \end{equation*}
Observe that $\dk^J_N  = \Span \bigl\{\t H,\, \t Z_{rs}: 1 \le r \le s \le N \bigr\}$, and set
\begin{equation*}
   \dm^J_N := \Span \bigl\{\t E,\, \t F,\, \t e_r,\, \t f_r: 1 \le r \le N \bigr\}.
\end{equation*}
Define $\tau: \dg^J_N \to \dg^J_N$ by $\tau(X) := \t X$.  We collect some of its properties:

\begin{lemma} \label{tilde auto} \begin{enumerate}
\item[(i)]  $\tau$ is an automorphism of\/ $\dg^J_N$ which preserves $\dsl_2$ and $\dr^J_N$.
\item[(ii)]  $\dg^J_N = \dk^J_N \oplus \dm^J_N$, and $\dm^J_N$ is the unique $\dk^J_N$-invariant complement of\/ $\dk^J_N$.
\item[(iii)]  The $\t H$-weights of\/ $\t E$, $\t F$, $\t e_r$, and $\t f_r$ are $2$, $-2$, $1$, and $-1$, respectively.
\item[(iv)]  The automorphism $\t\theta := \tau \circ \theta \circ \tau^{-1}$ of\/ $\dg^J_N$ acts by
\begin{equation*}
  \t\theta:\ \t H \mapsto -\t H, \quad \t E \mapsto -\t F \mapsto \t E, \quad
  \t e_r \mapsto \t f_r \mapsto -\t e_r, \quad \t Z_{rs} \mapsto \t Z_{rs}.
\end{equation*}
\item[(v)]  The homomorphism $\nu$ commutes with $\tau$, $\theta$, and $\t\theta$.  In particular,
\begin{eqnarray*}
   \nu(\t H) &:=& \t H - {\ts\frac{1}{4}} (\t f^T \t Z^{-1} \t e + \t e^T \t Z^{-1} \t f)
   \ =\ \t H - \oh (\t f^T \t Z^{-1} \t e - N), \\[6pt]
   \nu(\t E) &:=& \t E - {\ts\frac{1}{4}} \t e^T \t Z^{-1} \t e, \qquad
   \nu(\t F)\ :=\ \t F + {\ts\frac{1}{4}} \t f^T \t Z^{-1} \t f.
\end{eqnarray*}
\item[(vi)]  $\theta(\Omega_N) = \Omega_N$, $\theta(\Omega_{\dsl_2}) = \Omega_{\dsl_2}$, $\tau(\Omega_{\dsl_2}) = \Omega_{\dsl_2}$, and $\tau(\Omega_N) = \bigl( \frac{i}{2} \bigr)^N \Omega_N$.
\end{enumerate} \end{lemma}

\meno {\em Proof.\/}
Check~(v) directly for $\tau$ and $\theta$ on $E$ and $F$; the rest then follows from $[E,F] = H$ and the definition of $\t\theta$.  (If one proved an appropriate uniqueness property of $\nu$, (v)~would follow more elegantly.)  For~(vi), check that both $\theta$ and $\tau$ fix $\Omega_{\dsl_2}$ and then apply~(v).  We leave the rest to the reader.  $\Box$

\subsection{The slash actions $|_{k,L}$}

As described in the introduction, our results concern the \ido\ algebras of the $G^J_N$-\vb s over the homogeneous space $G^J_N/K^J_N = \bH \times \bC^N$.  These \vb s are in canonical bijection with \r s of $K^J_N$, and since the homogeneous space is simply connected, they are also in bijection with slash actions of $G^J_N$.  We shall consider only \irr\ $G^J_N$-\vb s.  Since $K^J_N$ is abelian, this means that we shall consider only line bundles, or in other words, scalar slash actions.

We will use the algebraic description of the IDO algebras developed by Helgason; see for example \cite{He77, He79, BCR12}.  This description depends only on the \r\ of $K^J_N$ defining the \vb, so we will not need the explicit descriptions of the slash actions given in \cite{CR}.  In fact, since $K^J_N$ is connected, we need only consider \r s of its Lie algebra $\dk^J_N$.  Helgason's algebraic machinery works for arbitrary \r s of $\dk^J_N$, regardless of whether or not they exponentiate to $K^J_N$; we will in fact treat all those with {\em invertible indices\/} (see below).  To put it differently, we will describe the IDO algebras of the scalar slash actions with invertible indices of all covers of $G^J_N$ on our homogeneous space; in particular, the metaplectic double cover.

The slash actions $|_{k,L}$ in \cite{CR} are defined for all $k \in \bZ$ and $L \in M_N^T(\bC)$.  In that paper the matrix $L$ is known as the {\em index\/} of the slash action, and in the number-theoretic applications it is taken to be positive definite symmetric with entries in $\oh\bZ$ and diagonal entries in $\bZ$; see Definition~2.9.  As discussed there in Section~5.4, the \r\ $\pi_{k,L}$ of $\dk^J_N$ associated to $|_{k,L}$ is 1-dimensional and acts by
\begin{equation*}
   \pi_{k,L}(\t H) := -k, \qquad \pi_{k,L}(\t Z_{rs}) := \pi L_{rs}.
\end{equation*}
We will write $\bC_{k,L}$ for the space of $\pi_{k,L}$.  The set of all \irrep s of $\dk^J_N$ is obtained by taking arbitrary $k \in \bC$ and $L \in M^T_N(\bC)$.  In the main results of this article we assume that $L$ is invertible.

\subsection{The IDO algebras}

Let us write $\bD_{k,L}$ for the algebra $\bD(G^J_N \times_{K^J_N} \bC_{k,L})$ of \dog s invariant with respect to the slash action $|_{k,L}$.  We now apply results stated in \cite{He77} to give an algebraic description of this algebra.  We will use the notation of Section~4 of \cite{BCR12}: by Theorem~4.1 of that paper,  $\bD_{k,L}$ is isomorphic to a certain algebra $\D_{\pi_{k,L}}$, which is defined by means of an auxiliary algebra $\E_{\pi_{k,L}}$.  For clarity we will suppress the $\pi$ in the subscript $\pi_{k,L}$.

Here $\pi_{k,L}$ is 1-dimensional, so $\E_{k,L}$ is simply $\dU(\dg^J_N)$.  The left ideal $\I_{k,L}$ of \cite{BCR12}~(34) is
\begin{equation*}
   \I_{k,L} := \dU(\dg^J_N) \bigl\{ Y + \pi_{k,L}(Y): Y \in \dk^J_N \bigr\},
\end{equation*}
the kernel of the canonical projection $\dU(\dg^J_N) \twoheadrightarrow \dU(\dg^J_N) \otimes_{\dk^J_N} \bC_{-k,-L}$.  It may be written in the more concrete form
\begin{equation} \label{IkL}
   \I_{k,L} := \dU(\dg^J_N) \bigl\{ \t H - k,\, \t Z_{rs} + \pi L_{rs}: 1 \le r \le s \le N \bigr\}.
\end{equation}

The situation simplifies because $K^J_N$ is connected and {\em reductive,\/} \ie\ $\dm^J_N$ is $\dk^J_N$-invariant: \cite{BCR12} Proposition~4.4 gives
\begin{equation} \label{DkL}
   \D_{k,L} = \bigl( \dU(\dg^J_N) \otimes_{\dk^J_N} \bC_{-k,-L} \bigr)^{\dk^J_N}
   = \dU(\dg^J_N)^{\dk^J_N} \big/ (\I_{k,L})^{\dk^J_N}
\end{equation}
(the superscript $\dk^J_N$ indicates the $\dk^J_N$-invariants).  Note that the algebra structure of $\D_{k,L}$ is unclear in the middle expression but clear in the right expression, because, as is easily checked, $(\I_{k,L})^{\dk^J_N}$ is a two-sided ideal in $\dU(\dg^J_N)^{\dk^J_N}$.

Henceforth we shall identify $\bD_{k,L}$ and $\D_{k,L}$.  Again by \cite{BCR12} Proposition~4.4, the graded commutative algebra of $\bD_{k,L}$ is $\S(\dm^J_N)^{\dk^J_N}$, the $\dk^J_N$-invariants of the symmetric algebra of $\dm^J_N$.  In light of Lemma~\ref{tilde auto}, this algebra has basis
\begin{equation} \label{basis1}
   \bigl\{ \t F^{i_F} \t E^{i_E} \t f^{I_f} \t e^{I_e}: i_F, i_E \in \bN;\,
   I_f, I_e \in \bN^N;\, 2i_F + |I_f| = 2i_E + |I_e| \bigr\}.
\end{equation}
Here $\t f^{I_f}$ and $\t e^{I_e}$ are given in multinomial notation:
\begin{equation*}
   I_f = (I_{f,1}, \ldots, I_{f,N}), \qquad
   \t f^{I_f} := \t f_1^{I_{f,1}} \cdots \t f_N^{I_{f,N}}, \qquad
   |I_f| := I_{f,1} + \cdots + I_{f,N}.
\end{equation*}

We will freely regard the monomials $\t F^{i_F} \t E^{i_E} \t f^{I_f} \t e^{I_e}$ as elements of $\S(\dm^J_N)$, $\dU(\dg^J_N)$, or $\bD_{k,L}$, depending on the context.  The next proposition describes multiplication in $\bD_{k,L}$; it follows from~(\ref{IkL}) and~(\ref{DkL}).  The subsequent results give conditions under which two IDO algebras are isomorphic.

\begin{prop} \label{DkL rules}
The IDO algebra $\bD_{k,L}$ has basis~(\ref{basis1}).  In order to write the $\bD_{k,L}$-product of two basis monomials as a linear combination of basis monomials, use the commutation relations in $\dU(\dg^J_N)$ with the following modifications:
\begin{itemize}
\item Replace each occurrence of $\t Z_{rs}$ by the scalar $-\pi L_{rs}$.
\item Move each occurrence of $\t H$ to the far right or far left according to the weight commutation rules in Lemma~\ref{tilde auto}~(iii), and then replace it by~$k$.
\end{itemize}
\end{prop}

\begin{lemma} \label{isos}
Suppose that $\mu$ is an automorphism of\/ $\dg^J_N$ such that $\mu(\dk^J_N) = \dk^J_N$ and $\pi_{k,L} = \pi_{k',L'} \circ \mu$.  Then $\mu(\dm^J_N) = \dm^J_N$ and $\bD_{k,L}$ is isomorphic to $\bD_{k',L'}$.
\end{lemma}

\meno {\em Proof.\/}
Since $\mu$ preserves $\dk^J_N$, it preserves $\dU(\dg^J_N)^{\dk^J_N}$.  Clearly it maps $\I_{k,L}$ to $\I_{k',L'}$, so by~(\ref{DkL}) it defines the desired isomorphism.  It preserves $\dm^J_N$ by Lemma~\ref{tilde auto}~(ii).  $\Box$

\begin{cor} \label{kL isos} \begin{enumerate}
\item[(i)]  $\bD_{k,L} \cong \bD_{-k,L}$ for all $k$ and $L$.
\item[(ii)]  $\bD_{k,L} \cong \bD_{k, MLM^T}$ for all $k$, $L$, and invertible $M$ in $M_N(\bC)$.
\item[(iii)]  $\bD_{k,L} \cong \bD_{k,L'}$ for all $k$ and all invertible $L$ and $L'$.
\end{enumerate} \end{cor}

\meno {\em Proof.\/}
For~(i), note that $\pi_{k,L} = \pi_{-k,L} \circ \t\theta$.  For~(ii), define $\mu_M: \dg^J_N \to \dg^J_N$ by
\begin{equation*}
   \mu_M|_{\dsl_2} := 1, \quad
   \mu_M(\t e_r) := (M \t e)_r, \quad \mu_M(\t f_r) := (M \t f)_r, \quad
   \mu_M(\t Z_{rs}) := (M \t Z M^T)_{rs}.
\end{equation*}
Check that $\mu_M$ is an automorphism of $\dg^J_N$ preserving $\dk^J_N$ such that $\pi_{k, MLM^T} = \pi_{k,L} \circ \mu_M$, so $\mu_M$ defines an isomorphism $\bD_{k, MLM^T} \to \bD_{k,L}$.  For~(c), apply~(b) and Takagi's factorization of complex symmetric matrices: since $L$ is invertible and symmetric, there exists $M$ such that $MLM^T = I_N$.  $\Box$

\medbreak
We conclude this section with an obvious but useful lemma which follows from the fact that since $Z_{rs} \mapsto 2\pi i L_{rs}$ in $\bD_{k,L}$, $\det(Z) \mapsto (2\pi i)^N \det(L)$.

\begin{lemma} \label{Z inverse}
For $L$ invertible, $\det(Z)^{-1} \mapsto \bigl(\frac{1}{2\pi i}\bigr)^N \det(L)^{-1}$ extends the projection $\dU(\dg^J_N)^{\dk^J_N} \twoheadrightarrow \bD_{k,L}$ to $\dU(\dg^J_N)^{\dk^J_N}[\det(Z)^{-1}]$.
\end{lemma}

\section{The centers $Z(\bD_{k,L})$ of the IDO algebras} \label{Z}

We now state our first main result: for all scalar slash actions with invertible indices, the center of the IDO algebra is polynomial in one variable and is generated by the {\em Casimir operator,\/} the image of the Casimir element $\Omega_N$.  As stated in the introduction, it was proven for $N=1$ in \cite{BCR12} and for $N=2$ in \cite{Da13}.

\begin{thm} \label{D center}
For $L$ invertible, the center $Z(\bD_{k,L})$ of\/ $\bD_{k,L}$ is the image of the center $\dZ(\dg^J_N)$ of\/ $\dU(\dg^J_N)$: the 1-variable polynomial algebra $\bC[\Omega_N + (\I_{k,L})^{\dk^J_N}]$.
\end{thm}

The proof of this theorem occupies the remainder of this section.  By Corollary~\ref{kL isos}, it suffices to prove it for any one invertible $L$.  The most convenient choice is $L = \frac{1}{2\pi} I_N$, so we define
\begin{equation*}
   \bD_k := \bD_{k,I_N/2\pi}, \qquad \I_k := \I_{k, I_N/2\pi}.
\end{equation*}

We remark that the Casimir operator of $\bD_{k,L}$ is denoted by $\C^{k,L}$ in~(2.4) of \cite{CR}.  Here it will be more convenient to use the equivalent operator
\begin{equation*}
   \C_{N,k} := \nu(\Omega_{\dsl_2}) + (\I_k)^{\dk^J_N} =
   (-i)^N \C^{k,I_N/2\pi} + {\ts\frac{1}{4}} N(N+4).
\end{equation*}

We know that the image of $\dZ(\dg^J_N)$ in $\bD_k$ is contained in $Z(\bD_k)$.  The image of the center $\dz(\dg^J_N)$ of $\dg^J_N$ is the scalars, so by Proposition~\ref{center} the image of $\dZ(\dg^J_N)$ is $\bC[\C_{N,k}]$.  Thus we come down to proving that $Z(\bD_k)$ is no bigger than $\bC[\C_{N,k}]$.

We begin by using Theorem~\ref{factorization} together with Lemma~\ref{tilde auto} to replace~(\ref{basis1}) by two more useful bases of $\bD_k$.  Let us write $\equiv$ for congruence modulo the ideal $(\I_k)^{\dk^J_N}$ defining $\bD_k$, and $\delta_{rs}$ for the Kronecker function.  For brevity, we set $X_\nu := \nu(X)$ for $X \in \dsl_2$.  Tracing definitions gives
\begin{equation} \label{formulas1}
   \t Z \equiv - \oh I_N, \quad
   \t H_\nu \equiv \t H + \t f^T \t e + \oh N, \quad
   \t E_\nu \equiv \t E + \oh \t e^T \t e, \quad
   \t F_\nu \equiv \t F - \oh \t f^T \t f.
\end{equation}

\begin{prop} \label{Dk rules1}
The IDO algebra $\bD_k$ has basis
\begin{equation} \label{basis2}
   \bigl\{ \t F_\nu^{i_F} \t E_\nu^{i_E} \t f^{I_f} \t e^{I_e}: i_F, i_E \in \bN;\,
   I_f, I_e \in \bN^N;\, 2i_F + |I_f| = 2i_E + |I_e| \bigr\}.
\end{equation}
In order to write the $\bD_k$-product of two basis monomials as a linear combination of basis monomials, use the following rules.  First, $X_\nu$ commutes with $\t f$ and $\t e$ for all $X \in \dsl_2$.  Second, $[\t E_\nu, \t F_\nu] = \t H_\nu$, and any occurrence of $\t H_\nu$ may be moved to the far right according to its weight commutation rules and replaced by $k + \t f^T \t e + \oh N$.  Third, $[\t e_r, \t f_s] \equiv \delta_{rs}$.  
\end{prop}

\meno {\em Proof.\/}
The fact that~(\ref{basis2}) is a basis follows from the fact that~(\ref{basis1}) is a basis: one may use an inductive argument based on the degree function on $\S(\dm^J_N)$ given by letting $\t F$ and $\t E$ have degree~2 and $\t f_r$ and $\t e_s$ have degree~1.  Theorem~\ref{factorization} shows that $\nu(\dsl_2)$ commutes with $\dr^J_N$.  The rest follows from~(\ref{basis2}).  $\Box$

\medbreak
Observe that by Lemma~\ref{tilde auto}, $\C_{N,k} \equiv \nu \circ \tau(\Omega_{\dsl_2}) =  \t H_\nu^2 + 2 \t H_\nu + 4 \t F_\nu \t E_\nu$.  The idea leading to the basis of $\bD_k$ given in the next proposition is that elements of~(\ref{basis2}) with both $i_F$ and $i_E$ non-zero can be replaced  with monomials involving $\C_{N,k}$ and at most one of $\t F_\nu$ and $\t E_\nu$.  The proof is left to the reader.

\begin{prop} \label{Dk rules2}
The IDO algebra $\bD_k$ has basis
\begin{equation*}
   \bigl\{ \t F_\nu^{i_F} \t E_\nu^{i_E} \t f^{I_f} \t e^{I_e} \C_{N,k}^{i_C}: i_F, i_E, i_C \in \bN;\,
   i_F i_E = 0;\, I_f, I_e \in \bN^N;\, 2i_F + |I_f| = 2i_E + |I_e| \bigr\}.
\end{equation*}
In order to write the $\bD_k$-product of two basis monomials as a linear combination of basis monomials, use the rules in Proposition~\ref{Dk rules1} together with the facts that $\C_{N,k}$ is central and $\t F_\nu \t E_\nu \equiv \frac{1}{4} (\C_{N,k} - \t H_\nu^2 - 2 \t H_\nu)$.
\end{prop}

As usual, write $(\bD_k)_{\Lie}$ for $\bD_k$ regarded as a Lie algebra via the commutator bracket.  The proof of Theorem~\ref{D center} hinges on a certain Lie subalgebra of $(\bD_k)_{\Lie}$.  The next lemma is proven by direct calculation.  Working in $\bD_k$, set
\begin{equation*}
   \E_{rs} :\equiv \oh (\t f_r \t e_s + \t e_s \t f_r), \quad
   \E := {\ts \sum_{r=1}^N} \E_{rr}, \quad
   \da_N := \Span \bigl\{ \E_{rs}: 1 \le r, s \le N \bigr\}.
\end{equation*}

\begin{lemma} \label{aN}
\begin{enumerate}
\item[(i)]  $\da_N$ is a Lie subalgebra of\/ $(\bD_k)_{\Lie}$.
\item[(ii)]  $\ep_{rs} \mapsto \E_{rs}$ is a Lie isomorphism from $\dgl_N$ to $\da_N$.
\item[(iii)]  $\E$ spans the center of\/ $\da_N$.
\item[(iv)]  The adjoint action $\ad(\E_{rs})$ of\/ $\E_{rs}$ on $\bD_k$ is the derivation determined by
\begin{equation*}
   \t F_\nu \mapsto 0, \quad \t E_\nu \mapsto 0, \quad
   \t f_i \mapsto \delta_{si} \t f_r, \quad \t e_i \mapsto -\delta_{ri} e_s, \quad
   \C_{N,k} \mapsto 0.
\end{equation*}
\item[(v)]  $\t F_\nu^{i_F} \t E_\nu^{i_E} \t f^{I_f} \t e^{I_e} \C_{N,k}^{i_C}$ is an eigenvector of\/ $\ad(\E)$ of eigenvalue $|I_f| - |I_e|$.
\end{enumerate} \end{lemma}

\meno{\em Proof of Theorem~\ref{D center}.\/}
First, the centralizer in $\bD_k$ of $\E$ is
\begin{equation*}
   (\bD_k)^\E = \Span \bigl\{ \t f^{I_f} \t e^{I_e} \C_{N,k}^{i_C}:
   I_f, I_e \in \bN^N;\, |I_f| = |I_e|;\, i_C \in \bN \bigr\}.
\end{equation*}
This is because by Lemma~\ref{aN}~(v), the only elements of the basis of Proposition~\ref{Dk rules2} commuting with $\E$ have $|I_f| = |I_e|$, which forces $i_F = i_E = 0$.

Second, the centralizer in $\bD_k$ of $\da_N$ is $(\bD_k)^{\da_N} = \bC[\C_{N,k}, \E]$.  To prove this, use Lemma~\ref{aN}~(ii) and~(iv) to check that the isomorphism from $\da_N$ to $\dgl_N$ carries the $\da_N$-module $\Span \{\t f^{I_f} \t e^{I_e}: |I_f| = |I_e| = m \}$ to the $\dgl_N$ module $\S^m(\bC^N) \ot \S^m(\bC^N)^*$, where $\S^m(\bC^N)$ is the $m^\thup$ symmetric power of the standard module of $\dgl_N$.  This is the tensor product of an irreducible module with its dual, so by Schur's lemma it contains a unique trivial submodule.  Equating dimensions gives
\begin{equation*}
   \bigl( \Span \bigl\{ \t f^{I_f} \t e^{I_e}: |I_f| = |I_e| \le m \bigr\} \bigr)^{\da_N} =
   \Span \bigl\{ \E^j: j \le m \bigr\}.
\end{equation*}

Finally, since $\C_{N,k}$ and $\E$ are clearly algebraically independent and $\D_k$ contains elements of non-zero $\E$-weight, we find that the center of $\bD_k$ is contained in, and therefore equal to, $\bC[\C_{N,k}]$.  $\Box$

\section{The characters of the IDO algebras in rank~$N>1$} \label{Characters}

Our second main result is the classification of the characters of $\bD_{k,L}$ for $L$ invertible.  By Corollary~\ref{kL isos} it suffices to treat $\bD_k$, just as in Section~\ref{Z}.  In this section we give the result for $N>1$.  We begin by defining two characters which exist for all $N$.  Let $\L_k^f$ be the subspace of $\bD_k$ spanned by elements ``with $\t F$'s and $\t f$'s on the left'', and let $\L_k^e$ be the subspace ``with $\t E$'s and $\t e$'s on the left'':
\begin{align*}
   \L_k^f :=& \Span \bigl\{ \t F_\nu^{i_F} \t f^{I_f} \t E_\nu^{i_E} \t e^{I_e}: 2i_F + |I_f| = 2i_E + |I_e| > 0 \bigr\}, \\[6pt]
   \L_k^e :=& \Span \bigl\{ \t E_\nu^{i_E} \t e^{I_e} \t F_\nu^{i_F} \t f^{I_f}: 2i_F + |I_f| = 2i_E + |I_e| > 0 \bigr\}.
\end{align*}

\begin{lemma} \label{LfLe}
\begin{enumerate}
\item[(i)]  $\L_k^f$ and $\L_k^e$ are two-sided ideals in $\bD_k$ of codimension~1.
\item[(ii)]  $\t\theta$ drops to an isomorphism $\t\theta_k: \bD_k \to \bD_{-k}$ such that $\t\theta_{-k} \circ \t\theta_k = 1$,
\begin{equation*}
   \t\theta_k(\t F_\nu^{i_F} \t f^{I_f} \t E_\nu^{i_E} \t e^{I_e}) =
   (-1)^{i_F + |I_f| + i_E}\, \t E_\nu^{i_F} \t e^{I_f} \t F_\nu^{i_E} \t f^{I_e}, \quad
   \t\theta_k(\L_k^f) = \L_{-k}^e.
\end{equation*}
\end{enumerate}
\end{lemma}

\meno {\em Proof.\/}
Note that $\L_f$ is the span of the elements of the basis~(\ref{basis2}) other than~1, because $\t f_r$ and $\t E_\nu$ commute.  Therefore it is of codimension~1.  By the multiplication rules given in Proposition~\ref{Dk rules1}, the product of two non-identity elements of~(\ref{basis2}) is a linear combination of non-identity basis elements.  Hence $\L_k^f$ is a two-sided ideal.

The existence of the isomorphism $\t\theta_k$ follows from Corollary~\ref{kL isos}~(i) and its proof.  The formula for its action on the elements of~(\ref{basis2}) follows from~(\ref{formulas1}) and Lemma~\ref{tilde auto}~(iv).  The fact that $|I_f| + |I_e|$ is even for all elements of~(\ref{basis2}) gives $\t\theta_{-k} \circ \t\theta_k = 1$.  It is now clear that $\t\theta_k(\L_k^f) = \L_{-k}^e$.  Part~(i) for $\L_k^e$ follows.  $\Box$

\medbreak
The next lemma gives generators for $\bD_k$.  It follows easily from Proposition~\ref{Dk rules2} and the fact that $\bD_k$ is commutative up to symbol.

\begin{lemma} \label{generators}
$\bD_k$ is generated by $\bigl\{ \t F_\nu \t e_r \t e_s,\, \t E_\nu \t f_r \t f_s,\, \E_{rs},\, \C_{N,k}: 1 \le r, s \le N \bigr\}$.
\end{lemma}

The two characters of $\bD_k$ existing for all $N$ are the projections to $\bC \cdot 1$ along $\L_k^f$ and $\L_k^e$.  We now give their actions on the generators in Lemma~\ref{generators}, and prove that for $N > 1$, they are the only characters of $\bD_k$.

\begin{prop} \label{2 chars}
The homomorphism $\chi_k^f: \bD_k \twoheadrightarrow \bC$ with kernel $\L_k^f$ maps
\begin{equation*}
   \t F_\nu \t e_r \t e_s \mapsto 0, \quad
   \t E_\nu \t f_r \t f_s \mapsto 0, \quad
   \E_{rs} \mapsto \oh \delta_{rs}, \quad
   \C_{N,k} \mapsto (k + \oh N) (k + \oh N +2).
\end{equation*}
The homomorphism $\chi_k^e: \bD_k \twoheadrightarrow \bC$ with kernel $\L_k^e$ is $\chi_{-k}^f \circ \t\theta_k$.  It maps
\begin{equation*}
   \t F_\nu \t e_r \t e_s \mapsto 0, \quad
   \t E_\nu \t f_r \t f_s \mapsto 0, \quad
   \E_{rs} \mapsto -\oh \delta_{rs}, \quad
   \C_{N,k} \mapsto (k - \oh N) (k - \oh N - 2).
\end{equation*}
\end{prop}

\meno {\em Proof.\/}
Since $\t F_\nu$ and $\t E_\nu$ commute with $\t f_r$ and $\t e_r$, both characters map $\t F_\nu \t e_r \t e_s$ and $\t E_\nu \t f_r \t f_s$ to zero.  The images of $\E_{rs}$ follow from
\begin{equation*}
   \E_{rs} \equiv \t f_r \t e_s + \oh \delta_{rs} \equiv \t e_s \t f_r - \oh \delta_{rs}.
\end{equation*}
For the images of $\C_{N,k}$, use $\t H \equiv k$ when on the far right, $\t H_\nu \equiv \t H + \E$, and
\begin{equation*}
   \C_{N,k} \equiv \t H_\nu (\t H_\nu + 2) + 4 \t F_\nu \t E_\nu = \t H_\nu (\t H_\nu - 2) + 4 \t E_\nu \t F_\nu.\ \ \Box
\end{equation*}

\begin{thm} \label{N chars}
For $N > 1$, $\chi_k^f$ and $\chi_k^e$ are the only characters of\/ $\bD_k$
\end{thm}

\meno {\em Proof.\/}
Let $\chi: \bD_k \to \bC$ be a character.  The kernel of $\chi$ must contain all commutators, and in particular all elements of non-zero $\E$-weight and all elements of the copy of $\dsl_N$ inside $\da_N$.  Therefore $\chi$ annihilates the generators $\t F_\nu \t e_r \t e_s$, $\t E_\nu \t f_r \t f_s$, and those $\E_{rs}$ with $r \not= s$, and maps all $\E_{rr}$ to the same value, say $x$.  Note that
\begin{equation*}
   \E_{12} \E_{21} \equiv \t f_1 \t e_2 \t f_2 \t e_1 \equiv \t f_1 \t e_1 \t e_2 \t f_2 \equiv (\E_{11} - \oh) (\E_{22} + \oh).
\end{equation*}
Since $\chi$ maps the left side to $0$ and the right side to $(x - \oh) (x + \oh)$, we have $x = \pm \oh$.

Suppose that $x = \oh$.  Note that
\begin{equation*}
   (\t F_\nu \t e_1 \t e_1) (\t E_\nu \t f_1 \t f_1)
   \equiv (\t F_\nu \t E_\nu) (\t e_1^2 \t f_1^2)
   \equiv (\t F_\nu \t E_\nu) (\t e_1 \t f_1) (\t e_1 \t f_1 + 1)
\end{equation*}
Since $\t e_1 \t f_1 \equiv \E_{11} + \oh$, $\chi$ maps the right side to $2 \t F_\nu \t E_\nu$.  It maps the left side to~$0$, so it annihilates $\t F_\nu \t E_\nu$.  Proceed as in the proof of Proposition~\ref{2 chars} to prove that $\chi$ matches $\chi_k^f$ on the generators in Lemma~\ref{generators}, so they are equal.

If $x = -\oh$, verify $\chi \circ \t\theta_{-k} = \chi_{-k}^f$ and deduce $\chi = \chi_k^e$, completing the proof.  $\Box$

\section{The representations of the IDO algebras in rank~$N=1$} \label{Reps}

In this section we study the case $N=1$: we classify the finite dimensional \r s of $\bD_k$.  In particular, we classify its characters: there are usually five, sometimes four.  We also prove that $\bD_k \cong \bD_{k'}$ \iff\ $k' = \pm k$.  Our strategy is to embed $\bD_k$ in the \uea\ $\dU(\dsl_2)$ and apply the techniques used to analyze \r s of $\dsl_2$.

Write $\t f$ for $\t f_1$, $\t e$ for $\t e_1$, and $\C_k$ for $\C_{1,k}$.  Recall that $\E \equiv \t f \t e + \oh$, and, when on the far right, $\t H \equiv k$.  Let us restate Proposition~\ref{Dk rules2} and Lemma~\ref{generators} in this context:

\begin{prop} \label{1generators}
For $N=1$, $\bD_k$ has generators $\bigl\{ \t F_\nu \t e^2,\, \t E_\nu \t f^2,\, \E,\, \C_k \bigr\}$, and basis  \begin{equation*}
   \bigl\{ \t F_\nu^{i_F} \t E_\nu^{i_E} \t f^{I_f} \t e^{I_e} \C_k^{i_C}:
   i_F, i_E, I_f, I_e, i_C \in \bN;\, i_F i_E = 0;\, 2i_F + |I_f| = 2i_E + |I_e| \bigr\}.
\end{equation*}
Here $\C_k$ is central, $\t f$ and $\t e$ commute with $\t F_\nu$ and $\t E_\nu$, $[\t e, \t f] \equiv 1$, and
\begin{equation*}
   [\t E_\nu, \t F_\nu] \equiv \t H_\nu \equiv \t H + \E, \qquad
   \t F_\nu \t E_\nu \equiv {\ts\frac{1}{4}} (\C_{N,k} - \t H_\nu^2 - 2 \t H_\nu).
\end{equation*}
\end{prop}

Keeping in mind the fact that $\bD_k$ is commutative up to symbol, we obtain:

\begin{cor} \label{1relations}
\begin{enumerate}
\item[(i)]  The following is a complete set of relations for the generators $\t F_\nu \t e^2$, $\t E_\nu \t f^2$, $\E$, and $\C_k$ of\/ $\bD_k$: $\C_k$ is central,
\begin{align}
   & \ad(\E):\, \t F_\nu \t e^2 \mapsto -2 \t F_\nu \t e^2, \quad
   \t E_\nu \t f^2 \mapsto 2 \t E_\nu \t f^2,
   \label{Ewgt} \\[6pt]
   & (\t F_\nu \t e^2) (\t E_\nu \t f^2) =
   {\ts\frac{1}{4}} (\E + \oh) (\E + {\ts\frac{3}{2}}) \bigl( \C_k - (\E + k) (\E + k + 2) \bigr),
   \label{FE} \\[6pt] \label{EF}
   & (\t E_\nu \t f^2) (\t F_\nu \t e^2) =
   {\ts\frac{1}{4}} (\E - \oh) (\E - {\ts\frac{3}{2}}) \bigl( \C_k - (\E + k) (\E + k - 2) \bigr).
\end{align}
\item[(ii)]  $\bD_k$ has basis $\bigl\{ (\t F_\nu \t e^2)^{i_F} (\t E_\nu \t f^2)^{i_E} \E^j \C_k^{i_C}:
   i_F, i_E, j, i_C \in \bN;\, i_F i_E = 0 \bigr\}$.
\smallbreak
\item[(iii)]  $\t\theta_k$ exchanges $\t F_\nu \t e^2$ and $-\t E_\nu \t f^2$ and maps $\E$ to $-\E$ and $\C_k$ to $\C_{-k}$.
\end{enumerate}
\end{cor}

We remark that applying $\t\theta_{-k}$ to~(\ref{FE}) in $\bD_{-k}$ gives~(\ref{EF}) in $\bD_k$.

At this point it is convenient to define an abstract copy of $\dsl_2$:
\begin{equation*}
   \dsl_2 := \Span \bigl\{ x, y, h \bigr\}; \qquad [x, y] = h, \quad [h, x] = 2x, \quad [h, y] = -2y.
\end{equation*}
Let $\omega$ be its Casimir element and, as in Section~\ref{Sec 2.1}, let $\theta$ be its Cartan involution:
\begin{equation*}
   \omega =  h^2 + 2h + 4yx; \qquad
   \theta: h \mapsto -h, \quad x \mapsto -y \mapsto x.
\end{equation*}

By Corollary~\ref{1relations}~(ii), we may define a linear map $\iota_k: \bD_k \to \dU(\dsl_2)$ by
\begin{equation*}
   \iota_k \bigl( (\t F_\nu \t e^2)^{i_F} (\t E_\nu \t f^2)^{i_E} \E^j \C_k^{i_C} \bigr) :=
   \bigl( y (h - k - \oh) \bigr)^{i_F} \bigl( x (h - k + \oh) \bigr)^{i_E} (h - k)^j \omega^{i_C}.
\end{equation*}

\begin{prop} \label{embedding}
$\iota_k$ is an embedding of\/ $\bD_k$ in $\dU(\dsl_2)$, and $\theta \circ \iota_k = \iota_{-k} \circ \t\theta_k$.
\end{prop}

\meno {\em Proof.\/}
$\iota_k$ is injective by the PBW theorem for $\dU(\dsl_2)$, so it suffices to check that it preserves the relations given in Corollary~\ref{1relations}~(i).  Clearly $\omega$ is central and $y (h - k - \oh)$ and $x (h - k + \oh)$ are of $\ad(h - k)$-weights $-2$ and $2$, respectively.  The image of~(\ref{FE}) holds because
\begin{align*}
   & \iota_k \bigl[ (\t F_\nu \t e^2) (\t E_\nu \t f^2) \bigr] = y (h - k - \oh) x (h - k + \oh) \\[6pt]
   &= {\ts\frac{1}{4}} (h - k + \oh) (h - k + {\ts\frac{3}{2}}) \bigl( \omega - h(h+2) \bigr) \\[6pt]
   &= \iota_k \bigl[ {\ts\frac{1}{4}} (\E + \oh) (\E + {\ts\frac{3}{2}}) \bigl( \C_k - (\E + k) (\E + k + 2) \bigr) \bigr].
\end{align*}

The reader may check $\theta \circ \iota_k = \iota_{-k} \circ \t\theta_k$.  Applying it to~(\ref{FE}) gives~(\ref{EF}).  $\Box$

\medbreak
Via the embedding $\iota_k$, any \r\ of $\dsl_2$ may be regarded as a \r\ of $\bD_k$.  For reference we give the images of the generators of $\bD_k$:
\begin{equation*}
   \iota_k: \t F_\nu \t e^2 \mapsto y (h - k - \oh), \quad
   \t E_\nu \t f^2 \mapsto x (h - k + \oh), \quad
   \E \mapsto h - k, \quad
   \C_k \mapsto \omega.
\end{equation*}
In keeping with the \r\ theory of $\dsl_2$, we make the following definition:

\meno {\bf Definition.}
A representation of $\bD_k$ is {\em admissible\/} if the action of $\E$ is semisimple with finite dimensional eigenspaces.

\subsection{Admissible representations} \label{Adm Reps}

We now proceed to classify the \irr\ admissible \r s of $\bD_k$, and to prove that they include its \irr\ finite dimensional \r s.  First we establish some notation.

Let $V$ be any \r\ of $\bD_k$.  We say that it has {\em central character\/} $c$ if $\C_k$ acts on it by the scalar $c$.  Because $\iota_k(\E + k) = h$, we define the {\em $\lambda$-weight space\/} $V_\lambda$ of $V$ to be the $\lambda$-eigenspace of $\E + k$.  We say that {\em $\lambda$ is a weight of\/} $V$ if $V_\lambda$ is non-zero, and we say that $V$ {\em has consecutive weights\/} if its weights are the intersection of a convex subset of $\bC$ with an additive coset of $2\bZ$.  On each such coset $\lambda + 2\bZ$ we put the obvious left-to-right order transferred from $2\bZ$: for any $\mu$ and $\mu'$ in $\lambda + 2\bZ$, we say $\mu > \mu'$ if $\mu - \lambda > \mu' - \lambda$.  For $c$, $\lambda$ in $\bC$, define
\begin{align*}
   M^-_k (c, \lambda) &:=
   \bigl\{ -\infty \bigr\} \cup \Bigl[ \bigl\{ k + \oh,\, k + {\ts\frac{3}{2}},\, 1 \pm \sqrt{c + 1} \bigr\}
   \cap \bigl( \lambda - 2\bN \bigr) \Bigr], \\[6pt]
   M^+_k(c, \lambda) &:=
   \bigl\{ \infty \bigr\} \cup \Bigl[ \bigl\{ k - \oh,\, k - {\ts\frac{3}{2}},\, -1 \pm \sqrt{c + 1} \bigr\}
   \cap \bigl( \lambda + 2\bN \bigr) \Bigr], \\[6pt]
   m^-_k (c, \lambda) &:= \max M^-_k(c, \lambda), \qquad
   m^+_k (c, \lambda) := \min M^+_k(c, \lambda), \\[6pt]
   \Delta_k(c, \lambda) &:= \bigl\{ \mu \in \lambda + 2\bZ: m^-_k(c, \lambda) \le \mu \le m^+_k(c, \lambda) \bigr\}.
\end{align*}

\begin{lemma} \label{Dk reps}
Suppose that $V$ is a \r\ of\/ $\bD_k$, $\lambda$ is a weight of\/ $V$, and $v_\lambda$ is an eigenvector of\/ $\C_k$ in $V_\lambda$ of eigenvalue $c$.  For $r \in \bN$, define
\begin{equation*}
   v_{\lambda - 2r} := (\t F_\nu \t e^2)^r v_\lambda, \qquad
   v_{\lambda + 2r} := (\t E_\nu \t f^2)^r v_\lambda.
\end{equation*}
\begin{enumerate}
\item[(i)]  The submodule\/ $\bD_k v_\lambda$ has central character $c$ and is\/ $\Span \{ v_{\lambda + 2s}: s \in \bZ \}$.
\item[(ii)]  $(\bD_k v_\lambda)_{\lambda + 2s} = \bC v_{\lambda + 2s}$: all weight spaces of\/ $\bD_k v_\lambda$ are 0- or 1-dimensional.
\item[(iii)]  $\bD_k v_\lambda$ has consecutive weights: there is an element\/ $m^-(v_\lambda)$ of\/ $M^-_k(c, \lambda)$ and an element\/ $m^+(v_\lambda)$ of\/ $M^+_k(c, \lambda)$ such that $v_\mu \not= 0$ \iff\
\begin{equation*}
   m^-(v_\lambda)\, \le\, \mu\, \le\, m^+(v_\lambda), \quad \mu \in \lambda + 2\bZ.
\end{equation*}
\item[(iv)]  The non-zero $v_\mu$ are a basis of\/ $\bD_k v_\lambda$.
\item[(v)]  Up to equivalence,\/ $\bD_k v_\lambda$ is determined by $\lambda$, $c$, $m^-(v_\lambda)$, and\/ $m^+(v_\lambda)$.
\item[(vi)]  $\bD_k v_\lambda$ has a unique \irr\ quotient, whose equivalence class is determined by $\lambda$ and $c$.  A basis of this quotient is given by the images of those $v_\mu$ such that $\mu \in \Delta_k(c, \lambda)$.
\end{enumerate}
\end{lemma}

\meno {\em Proof.\/}
For~(i), use Corollary~\ref{1relations}~(ii) and the facts that $\C_k$ is central and $v_\lambda$ is a joint eigenvector of $\E$ and $\C_k$.  For~(ii), use~(\ref{Ewgt}).

For~(iii) and~(iv), note that $v_{\lambda \pm 2r_0} = 0$ implies $v_{\lambda \pm 2r} = 0$ for $r \ge r_0$.  Hence there exist $m^-(v_\lambda)$ in $(\lambda - 2\bN) \cup \{-\infty\}$ and $m^+(v_\lambda)$ in $(\lambda + 2\bN) \cup \{\infty\}$ such that
\begin{equation} \label{basis of Dkvl}
   \bigl\{ v_\mu: \mu \in \lambda + 2\bZ,\, m^-(v_\lambda) \le \mu \le m^+(v_\lambda) \bigr\}
\end{equation}
is a basis of $\bD_k v_\lambda$.  Apply~(\ref{FE}) and~(\ref{EF}) to obtain
\begin{align}
   (\t F_\nu \t e^2) (\t E_\nu \t f^2) v_\mu &=
   {\ts\frac{1}{4}} (\mu - k + \oh) (\mu - k + {\ts\frac{3}{2}}) (c - \mu^2 - 2\mu) v_\mu,
   \label{lmax} \\[6pt] \label{lmin}
   (\t E_\nu \t f^2) (\t F_\nu \t e^2) v_{\mu} &=
   {\ts\frac{1}{4}} (\mu - k - \oh) (\mu - k - {\ts\frac{3}{2}}) (c - \mu^2 + 2\mu) v_\mu.
\end{align}
Suppose that $m^-(v_\lambda) > -\infty$, and write $\mu$ for it.  Then $v_{\mu - 2} = 0$, so~(\ref{lmin}) must be zero.  But $v_\mu \not= 0$, so either $\mu$ is $k + \oh$ or $k + \frac{3}{2}$, or $\mu^2 + 2\mu = c$.  Thus $m^-(c, \lambda) \in M^-_k(c, \lambda)$.  Similarly, use~(\ref{lmax}) to prove $m^+(v_\lambda) \in M^+_k(c, \lambda)$.

For~(v), since the basis~(\ref{basis of Dkvl}) is determined by $\lambda$, $m^-(v_\lambda)$, and $m^+(v_\lambda)$, it will suffice to prove that the action of $\bD_k$ on $v_\mu$ is determined by $\mu$ and $c$.  We know that $\C_k v_\mu = c v_\mu$ and $\E v_\mu = (\mu - k) v_\mu$.  By definition, $(\t F_\nu \t e^2) v_\mu = v_{\mu - 2}$ for $\mu \le \lambda$ and $(\t E_\nu \t f^2) v_\mu = v_{\mu + 2}$ for $\mu \ge \lambda$.  The action of $\t F_\nu e^2$ on those $v_\mu$ with $\mu > \lambda$ follows from~(\ref{lmax}), and the action of $\t E_\nu f^2$ on those $v_\mu$ with $\mu < \lambda$ follows from~(\ref{lmin}).

For~(vi), check that if $\mu \in \lambda + 2\bN$, then $\Span \{v_\nu: \nu \in \mu + \bZ^+ \}$ is $\bD_k$-invariant \iff~(\ref{lmax}) is zero, and that $m^+_k (c, \lambda)$ is the least such $\mu$.  Similarly, if $\mu \in \lambda - 2\bN$, then $\Span \{v_\nu: \nu \in \mu - \bZ^+ \}$ is $\bD_k$-invariant \iff~(\ref{lmin}) is zero, and $m^-_k (c, \lambda)$ is the greatest such $\mu$.  The result follows.  $\Box$

\begin{prop} \label{adm reps}
\begin{enumerate}
\item[(i)]  Every irreducible admissible \r\ of\/ $\bD_k$ has a central character, consecutive weights, and 1-dimensional weight spaces.
\item[(ii)]  Irreducible finite dimensional \r s of\/ $\bD_k$ are admissible.
\end{enumerate}
\end{prop}

\meno {\em Proof.\/}
If $V$ is a non-zero admissible \r, let $V_\lambda$ be a non-zero weight space.  Then $V_\lambda$ is finite dimensional and preserved by $\C_k$, so it contains an eigenvector $v_\lambda$.  If $V$ is \irr\ it must be equal to $\bD_k v_\lambda$.  Part~(i) now follows from Lemma~\ref{Dk reps}.  If $V$ is finite dimensional it contains eigenspaces of $\E$, so the same argument works.  $\Box$

\begin{thm} \label{1classification}
The central character and weights are complete invariants of the \irr\ admissible \r s of\/ $\bD_k$.  More precisely:
\begin{enumerate}
\item[(i)]  For all scalars $c$ and $\lambda$, there is an \irr\ admissible \r\ $V_k(c, \lambda)$ with central character $c$ and weights $\Delta_k(c, \lambda)$.
\item[(ii)]  $V_k(c, \lambda)$ and $V_k(c', \lambda')$ are equivalent \iff\ $c' = c$ and $\lambda' \in \Delta_k(c, \lambda)$.
\item[(iii)]  Every \irr\ admissible \r\ is equivalent to some $V_k(c, \lambda)$.
\end{enumerate}
\end{thm}

\meno {\em Proof.\/}
For~(i), let $V_k(c, \lambda)$ be a vector space with basis $\bigl\{v_\mu: \mu \in \Delta_k(c, \lambda) \bigr\}$.  Define an action of $\bD_k$ on $V_k(c,\lambda)$ as follows.  First, $\C_k v_\mu = c v_\mu$ and $\E v_\mu = (\mu - k) v_\mu$.  Second, $(\t E_\nu \t f^2) v_\mu$ is $v_{\mu + 2}$ for $\lambda \le \mu < m^+_k(c, \lambda)$, and $0$ for $\mu = m^+_k(c, \lambda)$, and similarly, $(\t F_\nu \t e^2) v_\mu$ is $v_{\mu - 2}$ for $\lambda \ge \mu > m^-_k(c, \lambda)$, and $0$ for $\mu = m^-_k(c, \lambda)$.  Third, use~(\ref{lmin}) to define $(\t E_\nu \t f^2) v_\mu$ for $\mu < \lambda$, and~(\ref{lmax}) to define $(\t F_\nu \t e^2) v_\mu$ for $\mu > \lambda$.

It is easy to check that this action satisfies the relations given in Corollary~\ref{1relations}, so $V_k(c, \lambda)$ is indeed a \r\ of $\bD_k$.  Its irreducibility together with Parts~(ii) and~(iii) of the theorem are immediate from Lemma~\ref{Dk reps}~(v) and~(vi).  $\Box$

\medbreak
As a corollary we obtain the characters of $\bD_k$.  For the proof, note that $V_k(c, \lambda)$ is 1-dimensional \iff\ $m^-_k(c, \lambda) = \lambda = m^+_k(c, \lambda)$.

\begin{cor} \label{1 chars}
$V_k(c, \lambda)$ is 1-dimensional \iff\ it is one of
\begin{equation*}
   V_k \bigl(0, 0\bigr), \quad
   V_k \bigl( (k \pm {\ts\frac{3}{2}})^2 - 1,\, k \pm \oh \bigl), \quad
   V_k \bigl( (k \pm {\ts\frac{5}{2}})^2 - 1,\, k \pm {\ts\frac{3}{2}} \bigl).
\end{equation*}
These five characters are distinct unless $k \in \bigl\{ \pm \oh,\, \pm \frac{3}{2} \bigr\}$, when two of them coincide.
\end{cor}

\medbreak
We remark that these four special values of $k$, being half-integral, correspond to genuine slash actions of the metaplectic group.  Also, for $N = 1$ the characters from Section~\ref{Characters} are
\begin{equation*}
   \chi_k^f = V_k \bigl( (k + {\ts\frac{3}{2}})^2 - 1,\, k + \oh \bigl), \quad
   \chi_k^e = V_k \bigl( (k - {\ts\frac{3}{2}})^2 - 1,\, k - \oh \bigl).
\end{equation*}

\subsection{Representations of $\dsl_2$}

Another way to prove the existence of $V_k(c, \lambda)$ is to realize it as a \sq\ of the restriction to $\iota_k(\bD_k)$ of an \iar\ of $\dsl_2$.  Let us briefly recapitulate $\dsl_2$-theory.  A \r\ $V$ of $\dsl_2$ is said to be {\em admissible\/} if $h$ acts semisimply with finite dimensional eigenspaces.  The eigenvalues of $h$ are called weights, and the $\lambda$-weight space is denoted by $V_\lambda$.  If $\omega$ acts by a scalar~$c$, $V$ is said to have central (or infinitesimal) character~$c$.

It is well-known that the \iar s of $\dsl_2$ have central characters, consecutive weights, and 1-dimensional weight spaces.  Moreover, they are classified by their central character and weights, and given any scalars $c$ and $\lambda$, there is up to equivalence a unique \iar\ with central character~$c$ having $\lambda$ as a weight.  Here is the classification:

\begin{itemize}
\item $L(\lambda)$ for $\lambda \in \bN$, the \irr\ finite dimensional \r\ with weights $\{ \lambda, \lambda - 2, \ldots, -\lambda \}$.  Its central character is $\lambda^2 + 2\lambda$.
\item $M^-(\lambda)$ for $\lambda \not\in \bN$, the \irr\ Verma module with weights $\lambda - 2\bN$.  Its central character is $\lambda^2 + 2\lambda$.
\item $M^+(\lambda)$ for $\lambda \not\in -\bN$, the \irr\ Verma module with weights $\lambda + 2\bN$.  Its central character is $\lambda^2 - 2\lambda$.
\item $P(c, \lambda)$ for $c + 1 \not= (\mu + 1)^2$ for any $\mu \in \lambda + 2\bZ$, the \irr\ principal series module with weights $\lambda + 2\bZ$ and central character $c$.
\end{itemize}

The following proposition gives the reductions of the restrictions to $\iota_k(\bD_k)$ of the \iar s of $\dsl_2$.  Its proof is an elementary application of the formula for $\iota_k$ and is omitted.

\begin{prop} \label{sl2 decomps}
Let $V$ be an \iar\ of $\dsl_2$ with central character $c$, weights $\Delta(V)$, and basis $\{v_\mu: \mu \in \Delta(V) \}$.  Regarded as an admissible \r\ of\/ $\bD_k$ via $\iota_k$, it remains \irr\ unless $\Delta(V)$ contains either $\{ k + \oh, k - \frac{3}{2} \}$ or $\{k - \oh, k + \frac{3}{2} \}$.
\begin{enumerate}
\item[(i)]  If $\{ k + \oh, k - \frac{3}{2} \} \subseteq \Delta(V)$, then $\{v_\mu: \mu \in \Delta(V), \mu \ge k + \oh \}$ is a\/ $\bD_k$-sub\r\ of\/ $V$ equivalent to $V_k(c, k + \oh)$.  The quotient of $V$ by this sub\r\ is $\bD_k$-equivalent to $V_k(c, k - \frac{3}{2})$.
\item[(ii)]  If $\{ k - \oh, k + \frac{3}{2} \} \subseteq \Delta(V)$, then $\{v_\mu: \mu \in \Delta(V), \mu \le k - \oh \}$ is a\/ $\bD_k$-sub\r\ of\/ $V$ equivalent to $V_k(c, k - \oh)$.  The quotient of $V$ by this sub\r\ is $\bD_k$-equivalent to $V_k(c, k + \frac{3}{2})$.
\end{enumerate}
\end{prop}

\subsection{Isomorphism classes}

Suppose that $\mu: \bD_k \to \bD_{k'}$ is an algebra isomorphism.  Since isomorphisms preserve centers, it must restrict to an isomorphism from $\bC[\C_k]$ to $\bC[\C_{k'}]$.  Therefore for some scalars $a \not= 0$ and $b$,
\begin{equation*}
   \mu(\C_k) = a \C_{k'} + b.
\end{equation*}

Isomorphisms also preserve adjoint-semisimplicity: if $\delta \in \bD_k$ and $\ad(\delta)$ acts semisimply on $\bD_k$, then $\ad\bigl( \mu(\delta) \bigr)$ acts semisimply on $\bD_{k'}$ with the same eigenvalues.  We claim that the space of ad-semisimple elements of $\bD_k$ is $\bC[\C_k] \oplus \bC \E$.  To prove this, check first that ad-semisimple elements must lie in the weight zero subalgebra $\bC[\C_k, \E]$.  Then consider the ad-action of elements in this subalgebra on the subspaces of $\bD_k$ of non-zero weight: for elements not in $\bC[\C_k] \oplus \bC \E$, it increases the total $(\C_k, \E)$-degree.  The claim follows.

Thus $\mu$ maps $\bC[\C_k] \oplus \bC \E$ to $\bC[\C_{k'}] \oplus \bC \E$.  Conclude that
\begin{equation*}
   \mu(\E) = d_0 \E + p(\C_{k'})
\end{equation*}
for some scalar $d_0 \not= 0$ and some polynomial $p(\C_{k'})$.  The set of eigenvalues of $\ad(\E)$ is $2\bZ$, so the same must be true of $\ad \bigl(d_0 \E + p(\C_{k'}) \bigr)$.  It follows that $d_0 = \pm 1$.  Observe that for $s \in \bZ$, $\mu$ sends the weight space $(\bD_k)_{2s}$ to $(\bD_{k'})_{2 d_0 s}$.

Recall from Lemma~\ref{LfLe}~(ii) and Corollary~\ref{1relations}~(iii) that there is an isomorphism $\t\theta_k: \bD_k \to \bD_{-k}$ sending $\E$ to $-\E$.  If $d_0 = -1$, then $\t\theta_{k'} \circ \mu$ is an isomorphism from $\bD_k$ to $\bD_{-k'}$ ``whose $d_0$ is $1$''.  Therefore, replacing $\mu$ by $\t\theta_k \circ \mu$ and $k'$ by $-k'$ if need be, we may assume that $d_0 = 1$.

Under this assumption $\mu$ sends the weight~2 space $(\bD_k)_2 = (\t E_\nu \t f^2) \bC[\C_k, \E]$ to $(\bD_{k'})_2$, so there are polynomials $q(\C_{k'}, \E)$ and $r(\C_k, \E)$ such that
\begin{equation*}
   \mu(\t E_\nu \t f^2) = (\t E_\nu \t f^2) q(\C_{k'}, \E), \qquad
   \mu^{-1}(\t E_\nu \t f^2) = (\t E_\nu \t f^2) r(\C_k, \E).
\end{equation*}
Applying $\mu$ to the second equation and distributing, we obtain
\begin{equation*}
   \t E_\nu \t f^2 = (\t E_\nu \t f^2) q(\C_{k'}, \E) r \bigl(a \C_{k'} + b, \E + p(\C_{k'})\bigr).
\end{equation*}
It follows that $q$ is a scalar $d_+$.  Similarly, $\mu(\t F_\nu \t e^2) = d_- \t F_\nu \t e^2$ for some scalar $d_-$.

Now apply $\mu$ to~(\ref{FE}) to obtain the following equation in $\bD_{k'}$:
\begin{equation*}
   d_- d_+ (\t F_\nu \t e^2)(\t E_\nu \t f^2) = {\ts\frac{1}{4}}
   (\E + p + \oh) (\E + p + {\ts\frac{3}{2}}) \bigl( (a \C_{k'} + b) - (\E + p + k) (\E + p + k + 2) \bigr).
\end{equation*}
Compute the product on the left using~(\ref{FE}) in $\bD_{k'}$ to arrive at
\begin{align*}
   & d_- d_+ (\E + \oh) (\E + {\ts\frac{3}{2}}) \bigl( \C_{k'} - (\E + k') (\E + k' + 2) \bigr) \\[6pt]
   = & (\E + p + \oh) (\E + p + {\ts\frac{3}{2}}) \bigl( (a \C_{k'} + b) - (\E + p + k) (\E + p + k + 2) \bigr).
\end{align*}
Comparing $\C_{k'}$-degrees gives $p(\C_{k'})$ scalar.  Comparing coefficients of $\E^2 \C_{k'}$ gives $d_- d_+ = a$, and then comparing coefficients of $\E \C_{k'}$ gives $p=0$.  Simplifying gives
\begin{equation*}
   a (\E + k') (\E + k' + 2) = (\E + k) (\E + k + 2) -b.
\end{equation*}
Equating all coefficients gives $a = 1$, $b = 0$, and $k = k'$.

Recall that for all $d \not= 0$, there is an automorphism $\mu_d$ of $\dsl_2$ defined by $x \mapsto d x$, $y \mapsto d^{-1} y$, and $h \mapsto h$.  Clearly $\mu_d$ preserves $\iota_k(\bD_k)$, so it may be viewed as an automorphism of $\bD_k$.  Let $\Aut(\bD_k)$ be the group of automorphisms of $\bD_k$.  We have proven the following result.

\begin{prop} \label{N1Dk isos}
\begin{enumerate}
\item[(i)]  $\bD_k \cong \bD_{k'}$ \iff\ $k = \pm k'$.
\item[(ii)]  For $k \not= 0$, $\Aut(\bD_k) = \{ \mu_d: d \in \bC^\times \}$, and $\mu_d \circ \mu_{d'} = \mu_{dd'}$.
\item[(iii)]  $\Aut(\bD_0) = \{ \mu_d,\, \t\theta_0 \circ \mu_d: d \in \bC^\times \}$, and $\t\theta_0 \circ \mu_d = \mu_{1/d} \circ \t\theta_0$.
\end{enumerate}
\end{prop}

\section{Remarks}  \label{Remarks}

We conclude with some directions for further investigation.  One possibility is the \r\ theory of $\bD_k$ for $N > 1$.  Comparing Sections~\ref{UrJN} and~\ref{Reps}, it would be natural to try to realize $\bD_k$ as a \sq\ of $\dU(\dsp_{2N})$, and thus obtain its \irrep s as a subset of those of $\dsp_{2N}$.  Note that $\bD_k$ cannot be realized as subalgebra of $\dU(\dsp_{2N})$, because it has no trivial \r.

Concerning the case of non-invertible indices, we have computed the characters of $\bD_{k,0}$ for $N=1$.  Probably this can be done for arbitrary $N$, but such indices are not of interest in number theory so we know of no reason to do so.

Let $G^J_1(\bC)$ be the {\em complex\/} Jacobi group of degree~1 and rank~1: the semidirect product of $\SL_2(\bC)$ acting on the 3-dimensional complex Heisenberg group.  The relevant stabilizer subgroup is $K^J_1(\bC) = \SU_2 \times \bC$, so up to equivalence the scalar slash actions are in bijection with the real characters of $\bC$.  In \cite{BCR12} it is proven that the center of the IDO algebra of each of these slash actions is the image of the center of $\dU \bigl(\dg^J_N(\bC) \bigr)$, and is generated by a pair of conjugate cubic Casimir operators.  It might be worthwhile to compute the \irrep s of these IDO algebras.  A deeper project would be to consider the \irr\ vector-valued slash actions corresponding to the higher dimensional \irrep s of $\SU_2$.

Perhaps the most interesting of the projects we mention here is the decomposition of the space of smooth sections $C^\infty_{\sec}(G^J_N \times_{K^J_N} \bC_{k, L})$ under the joint action of the commuting algebras $\dU(\dg^J_N)$ and $\bD_{k,L}$.  For $L$ invertible, both algebras commute with the Casimir operator $\C^{k,L}$, so they preserve its eigenspaces.  Consider the case $N = 1$: using Corollary~\ref{kL isos}~(iii) to move the results of Section~\ref{Reps} from $\bD_k$ to $\bD_{k,L}$, we see that generically, any joint eigenfunction of $\C^{k,L}$ and $\E = \oh (\t f \t e + \t e \t f)$ generates an \irrep\ of $\bD_{k,L}$.  Such a \r\ might generate an \irrep\ of $\dU(\dg^J_N) \otimes \bD_{k,L}$ under the action of $\dU(\dg^J_N)$.

Formulas for the IDOs $\C^{k,L}$ and $\E$ are given in \cite{BS98}: up to additive and multiplicative scalars, $\C^{k,L}$ is their Casimir operator (see pages~38 and~82), and $\E$ is their $\lambda(P_1)$ (see pages~59 and~61).  The operators $D_\pm$ and $\Delta_1$ they use to define the Casimir operator are essentially our $\t E_\nu$, $\t F_\nu$, and $\t H_\nu$.  Alternate formulas for $\C^{k,L}$ appear in \cite{Pi09}~(9), \cite{BR10}~(12), and \cite{CR}~(2.4) and~(2.6).

\def\eightit{\it} \def\bib{\bf} \bibliographystyle{amsalpha}

\end{document}